\documentclass[a4paper,12pt]{article}
\usepackage{amsfonts}
\usepackage{multirow}
\usepackage{graphicx}

\begin{document}

\title{Banishing divergence Part 2: Limits of oscillatory sequences, and applications.}
\author{David Alan Paterson\\
CSIRO CMSE\\
Graham Rd, Highett, 3090\\
Australia}
\date{\today}
\maketitle

\begin{abstract}
Sequences diverge either because they head off to infinity or because they oscillate. Part 1~\cite{Part1} of this paper laid the pure mathematics groundwork by defining Archimedean classes of infinite numbers as limits of smooth sequences. Part 2 follows that with applied mathematics, showing that general sequences can usually be converted into smooth sequences, and thus have a well-defined limit.

Each general sequence is split into the sum of smooth, periodic (including Lebesgue integrable), chaotic and random components. The mean of each of these components divided by a smooth sequence, or the mean of the mean, will usually be a smooth sequence, and so the oscillatory sequence will have at least a leading term limit. Examples illustrate the wide range of oscillatory sequences that have leading term limits. 

Methologies are given for applications of such limits in four other contexts. One is for finding the limit of a misbehaving function at a point on the real number line. A second uses a nonstandard type of contour integration to find the limit of a function on the complex numbers. The third application is to Riemann sums for evaluating improper integrals. The final application is in evaluating the square of the Dirac delta function.
\end{abstract}
\section{Introduction}

This paper is the first step towards a mathematical formalism in which every infinite sequence of real numbers has a limit. Part 1~\cite{Part1} defined infinite and infinitesimal numbers in a new way, as limits of a specific "smooth" subset of sequences.

Part 2 handles handles limits of general sequences and in particular oscillatory sequences $\lim_{n\to\omega}f(n)$. Each sequence is split into the sum of smooth $s$, periodic $p$ (including Lebesgue integrable), chaotic $k$ and random $r$ components as follows: $f(n)={}^1s(n)+\sum_j{}^js(n){}^jp(n)+\sum_k{}^ks(n)$ ${}^kk(n)+\sum_l{}^lr(n)$ where each appearance of $s$ can be a different smooth sequence. The mean of each of these components, or the mean of the mean, will usually be a smooth sequence. Care is taken to ensure that each mean is definied locally at a sufficiently large $n$ rather than being smeared over several values of $n$.

This use of the mean maps general sequences onto smooth sequences and, since smooth sequences have a limit by Part 1 of this paper, many if not all sequences have a limit. Examples of difficult sequences with well-defined limits are presented. 

This paper then presents methologies for applications of such limits on the real $x, a$ and complex $z, z_0$ numbers $\lim_{x\to a}f(x)$, $\lim_{z\to z_0}f(z)$, for an application to improper integrals as Riemann sums $\int_0^{\omega}f(x)dx=\lim_{n\to\omega^2}$ $\frac{1}{\sqrt{n}}\sum_{k=1}^nf((k-.5)/\sqrt{n})$ and similar, and the paper finishes with a brief final note about squaring the Dirac delta function.

All of the work presented in this paper requires only elementary mathematics. All the techniques would have been accessible to any mathematician living 100 years ago, and indeed five separate pivotal papers on which Part 1 rests were published between 1895 and 1907. If this work hasn't appeared in print before, then the only reason I can think of is that this paper touches very briefly on a remarkably broad range of mathematical techniques, including infinite ordinal numbers, big O notation, axioms of set theory, Box-Jenkins time series analysis, probability and the ensemble mean, strange attractors, real and complex analysis and generalized functions.

As much as possible of the work presented here is new. 

\section{Notation from Part 1}

The following is a summary of the notation defined in Part 1~\cite{Part1} that is needed by Part 2.

\emph{Finite commutativity} resolves many classical paradoxes involving infinity by rejecting the process on bijection~\cite{bijection} on ordered sets.

$\omega$ is the first infinite ordinal number. What matters for this paper is that $\omega$ is bigger than all integers and that it's defined to be commutative: $\omega \neq 1 + \omega = \omega + 1$. A typical infinitesimal number is $1/\omega$.

An infinite \emph{sequence} of real numbers is a map $s$ from the natural numbers to the real numbers. Write as $s(n)$ or $s_n$.

A \emph{prototype} is a special type of sequence, denoted $p_n$ or $p$. No two prototypes have the same asymptotic behavior at infinity. These prototypes include such sequences as:
$1/n^3, n, n^\alpha, \ln(n), \exp(n), n^\alpha \ln(n), n^\alpha/\ln(\ln(n)),$ $\exp(\alpha n^{0.7})$ etc. for each nonzero real number $\alpha$. Many prototypes will already be familiar to those who use Big O notation~\cite{BigO}. The real numbers are based on the prototypes 1 and 0.

A \emph{leading term limit} is a form of asymptotic limit. The definition is: If for every positive real number $\epsilon$, there is an integer $N_\epsilon$, a non-zero real number $c$, and a prototype $p$, such that $|(s_n/p_n)-c|<\epsilon$ for all integer $n\ge N_\epsilon$, then it is said that the leading term limit $^1\lim s=cp$. If, instead, $s_n=0$ for all $n \ge N_\epsilon$ then $^1\lim s=0$.

The \emph{second term limit} is the leading term limit of $s_n-cp_n$ and is written $^2\lim s$.

A \emph{limit} is a finite sum of leading, second and higher term limits and is written $\lim s=\sum {}^ic{}^ip$.

A \emph{Cauchy limit} is the ordinary everyday definition of limit. If the Cauchy limit exists then it differs from the leading term limit only when the leading term limit is an infinitesimal.

A \emph{leading limitable} sequence is a sequence with a leading term limit according to the definition above. Sequences that are not leading limitable may still have a leading term limit, but only if they can be equivalenced to a leading limitable sequence.

A \emph{smooth} sequence is a sequence with a limit according to the definition above. Sequences that are not smooth may still have a limit, but only if they can be equivalenced to a smooth sequence.

$\mathbb{I}_n$ is used to denote an ordered field defined by a ratio of a pair of smooth sequences. It contains only infinite and infinitesimal and real numbers.

$^*\mathbb{R}$ is the field of hyperreals from non-standard analysis~\cite{2} that includes both infinite and infinitesimal numbers.

\section{General sequences}

Any sequence that is not smooth is called here \emph{oscillatory}. A general guideline, though by no means a precise mathematical statement, is that the limits of smooth sequences can be written:
\begin{displaymath}
\lim_{x\to\omega}f(x)=f(\omega)
\end{displaymath}
whereas limits of oscillatory sequences can't.

Care has to be taken because in the literature there are examples like~\cite{2, 25}:
\begin{displaymath}
-1/12=\sum_{n=1}^\infty n \mbox{ ;   } -1=\sum_{n=0}^\infty 2^n \mbox{ and }\frac{2}{3}=\sum_{n=0}^\infty (-1)^n
\end{displaymath}
The method used here does not yield counterintuitive results like these. Because of ambiguities in the treatment of infinity using this notation, largely because it allows inappropriate use of bijection, this notation won't be used in this paper again.\\

The following method for decomposing sequences into components was inspired by Box-Jenkins time series analysis~\cite{26}, but is made much easier because here there is a predefined map whereas in time series analysis the map is unknown and has to be determined. In Box-Jenkins analysis a sequence of measured data will typically be split into a trend, a periodic component and a random component. The periodic and random components themselves will in the most general case be non-stationary, with a varying amplitude.

In the general sequence $f_n$, start by taking $n$ to be an integer that is sufficiently large for the overall trend to be smooth, and for the oscillatory component to have an amplitude that varies smoothly. (This discription is deliberately vague, it will be firmed up after presenting examples). Only after this general sequence is mapped onto a smooth sequence will the limit $n\to\omega$ be equivalenced to that of the smooth sequence. The sequence could be written as:
\begin{displaymath}
f_n={}^1s_n+{}^2s_nb_n+r_n
\end{displaymath}
where $^1s$ and $^2s$ are smooth mappings, $b$ is a deterministic oscillatory mapping and $r$ is random. I would like to use the word ``ergodic"~\cite{27,ensemble} to describe $b$, but that's too restrictive, as will be seen in the examples below. If this sequence can be decomposed as above, and if it is possible to calculate the mean value of $b_n$ by the integral over the period or similar methods, and the ``ensemble mean"~\cite{ensemble} value of $r_n$ by probabilistic methods, then
\begin{displaymath}
\bar{f}_n={}^1s_n+{}^2s_n\bar{p}_n+\bar{r}_n
\end{displaymath}
and the resulting mean of $f$ will be smooth enough for the existence of a leading term limit, at least. The decomposition won't necessarily be unique, but that doesn't matter because the definition of limit on a smooth sequence ensures that it is closed under the operations of term by term addition and multiplication.\\

This introduction to oscillatory sequences suffices for the proof of an important lemma.

\emph{Lemma} $\mathbb{I}_n$ is a proper subset of $^*\mathbb{R}$

\emph{Proof} $^*\mathbb{R}$ is known to be the largest field that can be constructed from sequences~\cite{2} and  $\mathbb{I}_n$ is an ordered field constructed from sequences~\cite{Part1} so $\mathbb{I}_n \subset{}^*\mathbb{R}$. Let $s_n=0,1,0,1,0,1,\ldots$ and $t_n=1,0,1,0,1,0,\ldots$ and $u_n=.5,.5,.5,.5,.5,.5,\ldots$. Then under the action of ultrafilter $\mathfrak{U}$ in $^*\mathbb{R}$ we have $s=_{\mathfrak{U}} t\ne _{\mathfrak{U}} u$ but under the action of taking the mean of an oscillatory sequence in $\mathbb{I}_n$, $\lim s=\lim t=\lim u$ so  $\mathbb{I}_n\ne{}^*\mathbb{R}$.\\

Here are some examples of limits on difficult sequences. In each case the aim is to find a leading term limit, because the sequence minus the leading term limit yields a second sequence that can be used to find a second term limit, etc. Examples 1 and 2 were chosen specifically as cases where Ces\`{a}ro summation~\cite{28}  fails.\\

\emph{Example} 1.

$b_ns_n=(-\mathrm{e})^n$.

Write this as $b_n=(-1)^n$, $s_n=\mathrm{e}^n$ and note that $(-1)^n$ is bounded and oscillatory. $n$ has a 50\% chance of being even and 50\% chance of being odd so the mean of the bounded oscillatory component is $0.5\times(1)+0.5\times(-1)=0$ and $\bar{b}_ns_n=0\times e^n=0$.

\emph{Example} 2.

$b_n=(-1)^{\lfloor\ln (n)\rfloor}$.

For sufficiently large $n$, the floor of $\ln(n)$ has a 50\% chance of being even and 50\% chance of being odd so $\bar{b}_n=0.5\times(1)+0.5\times(-1)=0$.

\emph{Example} 3.
 
$b_n=0,1,0,0,1,0,0,0,1,0,0,0,0,1,0,\ldots$

The probability of $b_n$ having the value 1 can be calculated and turns out to be very close to $1/(\sqrt{2n}-.5)$. So $\bar{b}_n=1/(\sqrt{2n}-.5)$ and the leading term limit is $1/\sqrt{2n}$.

\emph{Example} 4.

$b_n=\sum_{j=1}^{m}\sin(d_jn+e_j)$ where all $d_j$ are nonzero.

For sufficiently large $n$, each of the sinusoidal terms can be integrated over its period to get a mean which in this case happens to be zero, so $\bar{b}_n=\sum_{j=1}^{m}0=0$.

\emph{Example} 5.

$f_n$ is an element of the normal distribution $N(0, s(n))$.

An element of $N(0, s(n))$ is the same as $s(n)$ times an element of $N(0,1)$ which is oscillatory with mean zero, so $\bar{f}_n=s(n)\times 0=0$.

\emph{Example} 6.

$b_n$ is the $x$ or $y$ component of the H\'{e}non strange-attractor~\cite{29} with $x_{i+1}=y_i+1-1.4x_i^2$ and $y_{i+1}=0.3x_i$.

This has a well-defined mean with $\bar{x}_n=0.2573\ldots$ and $\bar{y}_n=0.07719\ldots$

\emph{Example} 7.

Use $n$ to define a random real number $r_n \in U[0, 1)$ and set $b_n=0$ when $r_n$ is rational and $b_n=1$ when $r_n$ is irrational.

Using Lebesgue integration~\cite{lebesgue} the mean is found to be one, so $\bar{b}_n=1$.

\emph{Example} 8.

$f_n$ is the number of primes less than or equal to $n$.

Write $f_n$ as $s_n(f_n/\mathrm{Li}(n))$ where $s_n = \mathrm{Li}(n)$ is the logarithmic integral function. $b_n=f_n/\mathrm{Li}(n)$ is bounded and known to have a mean value of one~\cite{30}. The leading term limit comes from $\mathrm{Li}(n) \sim \ln(n)/n$.

\emph{Example} 9.

$f_n=\tan(n)=b_ns_n$. This is tough because not only is tan unbounded, its standard deviation is also unbounded and the expected value of the standard deviation actually increases with $n$. Set $s_n=\sqrt{n}$ and $b_n=\tan(n)/\sqrt{n}$. Then $b_n$ is oscillatory with a bounded oscillatory standard deviation that does not tend to zero. By symmetry of $\tan(n)$ about zero the mean $\bar{b}_n=0$.

\emph{Example} 10.
Two closely related examples with the same leading term limit.

10a. $f_n = \tan({n_{\pmod{\pi/2}}})$

10b. $x \in U[0, \pi/2)$ ,  $f_n = \tan(x)$

This was sufficiently difficult that it was worth checking three ways, by direct and Monte-Carlo simulation, by integration over a period, and by use of order statistics~\cite{orderstats}. Despite the function, median, quartiles etc. being time independent, the mean is time dependent because with an increasing number of values of $n$ the singularity at $\pi/2$ is approached more and more closely.
\begin{figure}[ht]
\centering
\includegraphics{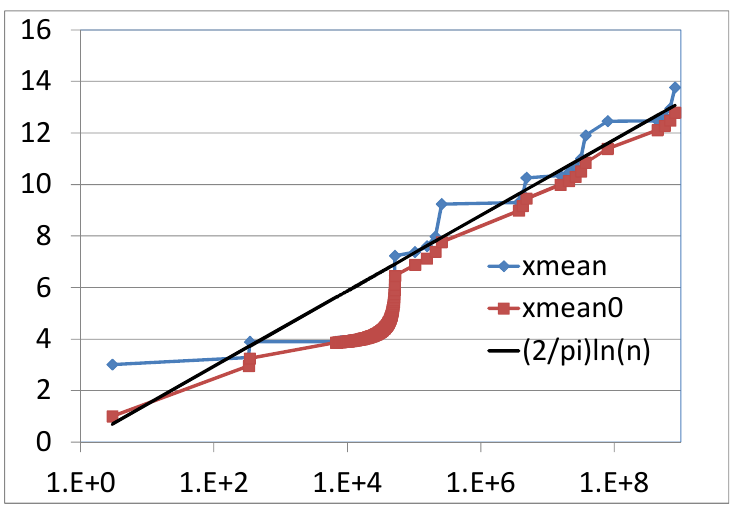}
\caption{Direct simulation of 10a. Upper and lower values of the mean.}
\end{figure}

Order statistics of 10b yields $\frac{1}{n}\sum_{k=1}^n\tan(k\pi/2(n+1))$.

Integration over a period for both 10a and 10b yields $\int_0^{1-\alpha/n}\tan(2x/\pi)dx$ where $\alpha \approx 1$.

All three methods give leading term limit $^1\lim f(i)=(2/\pi)\ln(n)$.\\

From the above ten examples it can be seen that the limit of the oscillatory component $b_n$ as defined above at sufficiently large but finite number $n$ exists for a remarkably large range of types of oscillatory sequences. In every example above $\bar{f}_n$ is a smooth sequence.\\

\emph{Conjecture}: Every sequence has a well-defined leading limit.

The above examples provide enough information to start attacking this conjecture. The following is a generalized methodology for generating such a limit.
\begin{itemize}
\item	Recall that $f_n = f(n)$ where $f$ is a mapping.
\end{itemize}
\emph{Periodics}
\begin{itemize}
\item	Start with simply periodic functions where there exists a real number $x_0$ such that $f(x+x_0)=f(x)$ for all necessary $x$.
\item	Fold the integers at period $x_0$ to get set $F\subset[0,x_0)$ and carry out a Lebesgue integration~\cite{lebesgue} of $f(x)$ over $F$.
\item	If the Lebesgue integral exists then the mean $\bar{f}$ is the Lebesgue integral divided by the Lebesgue measure of $F$.
\item	If the Lebesgue integral does not exist because of singularities in $f(x)$ at a finite number of points in $F\subset[0,x_0)$ then use order statistics~\cite{orderstats} to calculate how close $x$ gets to those singularities as a function of $n$ and use that expected distance as a bound on the integral. The integral will then have a limit on $\mathbb{I}_n$ as a function of $n$ and so the mean $\bar{f}(n)$ defined as the Lebesgue integral over the Lebesgue measure of $F$ will be well-defined.
\item	The next class of sequences comes from transformed periodic mappings $f(g(x))$ where $f$ is periodic and $g$ is leading term limitable. The mean of $f(g(x))$ is taken to be equal to the mean of $f(x)$.
\item	The next class of sequences $h(n)$ is of those that tend to a simply periodic or transformed periodic mapping. If $f(g(n))$ has a mean that is leading term limitable using prototype $p$ then this is equivalent to saying that the Cauchy limit $\lim_{n\to\infty}((h(n)-f(g(n))/p_n)=0$. The mean of $h(n)$ is then taken be the mean of the simply periodic or transformed periodic sequence that it tends to. The phrase ``periodic sequence" below means $h(n)$.
\end{itemize}
\emph{Randoms}
\begin{itemize}
\item	Let $f(i)$ be random. The mean and other statistical parameters may vary arbitrarily with $i$. From a sufficiently large number of actualisations take the ensemble mean~\cite{ensemble} at given $i = n$.
\item	If the ensemble mean is not easy to find (which sometimes happens when $f(i)$ is unbounded) then use order statistics to estimate how large $f(i)$ is likely to be at given $i$. This will likely allow an integral to be found for the ensemble mean.
\item	Treat transformed randoms and sequences that tend to randoms in the same way as for periodics.
\end{itemize}
\emph{Chaos}
\begin{itemize}
\item	Chaotic sequences will tend to a strange attractor. The best-known strange attractors have a well-defined mean.
\item	Treat variations on this in exactly the same way as variations on periodic functions.
\end{itemize}
\emph{Sums}
\begin{itemize}
\item	The method used to calculate the ensemble mean of a random sequence at given $n$ will probably generate a smooth sequence. The product of two smooth sequences is a smooth sequence. For full generality of oscillatory sequences each periodic sequence and chaotic sequence needs to be multiplied by a smooth sequence to recover smooth trends in amplitude. The limit of the product of smooth sequences is the product of the limits.
\item	The ensemble mean of a random sequence as a function of $n$ is not necessarily smooth. If a non-smooth ensemble mean can be converted to a finite sum of smooth times periodic and smooth times chaotic sequences (term by term multiplication and addition) then a well-defined limit exists. If not then some other method would be needed to find the limit, such as the imposition of symmetry conditions. By definition, the ensemble mean of an ensemble mean is identical to the original ensemble mean.
\item	So to summarise, it is very likely that a general sequence has a well-defined limit if it can be decomposed into a finite sum of: smooth sequence, smooth sequence times periodic sequence, smooth sequence times chaotic sequence, and random sequence -- providing the ensemble mean of the random sequence can be decomposed into a finite sum of: smooth sequence, smooth sequence times periodic sequence, and smooth sequence times chaotic sequence.
\end{itemize}

\emph{Lemma}
The complete set of sequences that can be mapped by this method to a smooth sequence are closed under addition and multiplication.

\emph{Proof}
Under addition because the sum of two finite sums is a finite sum. Under multiplication because, temporarily leaving aside the issue of inverses: the product of two smooth sequences is smooth, the product of two periodic sequences is periodic, the product of two chaotic sequences is chaotic, the product of two random sequences is random, the product of a periodic and chaotic sequence is chaotic, the product of a smooth and random sequence is random with smooth ensemble mean, the product of a periodic and random sequence is random with periodic ensemble mean, and the product of a chaotic and random sequence is random with chaotic ensemble mean, and the cases of smooth times periodic and smooth times chaotic are explicitly taken care of in the definition. When inverses are present, products become simplified.\\

In mathematical notation:\\

a) Split the sequence $f$ into smooth $s$, periodic $p$, chaotic $k$ and random $r$ components:
\begin{displaymath}
f(n)={}^1s(n)+\sum_j{}^js(n){}^jp(n)+\sum_k{}^ks(n){}^kk(n)+\sum_l{}^lr(n)
\end{displaymath}
where each sum is finite and each appearance of $s$ can mean a different smooth sequence.\\

b) Take the ensemble mean of each random component and split it into smooth, periodic and chaotic components:
\begin{displaymath}
^l\bar{r}(n)={}^ls(n)+\sum_m{}^ms(n){}^mp(n)+\sum_i{}^is(n){}^ik(n)
\end{displaymath}
where again each sum is finite and each appearance of a sequence can be a different sequence.\\

c) Take the mean of each periodic and chaotic component. The mean will be dependent on $n$ if and only if order statistics are required. Let's write it as if it is independent of $n$.
\begin{eqnarray*}
&^l\bar{\bar{r}}(n)={}^ls(n)+\sum_m{}^ms(n){}^m\bar{p}+\sum_i{}^is(n){}^i\bar{k} \\
&\bar{f}(n)={}^1s(n)+\sum_j{}^js(n){}^j\bar{p}+\sum_k{}^ks(n){}^k\bar{k}+\sum_l{}^l\bar{\bar{r}}(n)
\end{eqnarray*}
 
 d) Then $\bar{f}(n)$ will be a finite sum of smooth functions which is usually smooth, and we can take the limit using the techniques described earlier:
\begin{displaymath}
\lim_{n\to\omega}f(n)=\lim_{n\to\omega}\bar{f}(n)
\end{displaymath}
This gives the mapping from a general sequence to a smooth sequence which allows euivalencing of the limits of the two.

Note that periodic and chaotic components are treated in exactly the same way, except that the mean of a strange attractor is not calculated from the integral over a period. For practical use, the leading term limits of each component can be used to get the leading term limit of $f$ and then that subtracted from $f$ gives a new sequence that can be used to give the second term limit, etc.

\section{Limits on the real and complex numbers}

What about the limit of a function at a real number?

When a function $f$ is continuous at $x = a \in \mathbb{R}$ then there's no problem. $\lim_{x\to a}f(x)=f(a)$.

When $f$ is not continuous at $x = a$ then the following method works:\\

Allow $\epsilon$ to tend to a specific infinitesimal $\epsilon_0$. Split $f(x)$ near $a$ into $f(x)=f_1(x)+f_2(x)f_3(x)$ where $f_1$ and $f_2$ are monotone continuous near $a + \epsilon$ and/or $a - \epsilon$ and $f_3$ is ergodic~\cite{27} at the same point(s). Take the mean value of the ergodic component to get $\bar{f}_3$. If $f(x)$ has a value on the real number line at both $a + \epsilon$ and $a - \epsilon$  then use:

$\lim_{x\to a}f(x)=\lim_{\epsilon\to\epsilon_0}(1/2)(f_1(a+\epsilon)+f_2(a+\epsilon)\bar{f}_3(a+\epsilon)+f_1(a-\epsilon)+f(2(a-\epsilon)\bar{f}_3(a-\epsilon))$

Otherwise, if $f$ doesn't exist on the real numbers at either $a + \epsilon$ or $a - \epsilon$, use the one-sided limit:

$\lim_{x\to a}f(x)=\lim_{\epsilon\to\epsilon_0}(f_1(a\pm\epsilon)+f_2(a\pm\epsilon)\bar{f}_3(a\pm\epsilon)$\\

This leads to results such as the following.
\begin{eqnarray*}
&\lim_{x\to\pi/2}\tan(x)=0\qquad\lim_{x\to\epsilon_0}\ln(x)=-\ln(1/\epsilon_0)\\ &^1\lim_{x\to 1/\omega}|1/\sin(x)|=\omega\qquad\lim_{x\to 0}\cos(1/x)/x=0
\end{eqnarray*}
What about the complex numbers? There's no trouble with writing a complex number as $z=a+ib$ where $a$ and $b$ are infinite or infinitesimal numbers, but this doesn't immediately help much in determining the infinite limits of useful complex functions.\\

Let's look at the infinite limits of $e^z$. This is smooth and the limits are $e^{\omega}$ and $e^{-\omega}$ at the ends of the real number line. $e^z$ is periodic with mean zero in the imaginary direction so has limit 0 at both ends of the imaginary number line. Further, $e^z$ is periodic with mean zero along every straight line passing through the origin other than the real number line.

To get the transition from $e^{\omega}$ to 0 consider the complex number $z=\omega+ib$ for real b.

$e^{\omega+ib}=e^{\omega}(\cos(b)+i\sin(b))$.

The mean of a function is the integral (along a contour at constant distance from $z=0$) divided by the interval. So the limit of $e^{\omega+ib}$ is

$e^{\omega}(\sin(b)+i(1-\cos(b)))/b$.\\

What about the limit of a complex function at a finite complex number? On the complex numbers, calculate the mean from the average around a contour at infinitesimal radius $\epsilon_0$.

$\lim_{z\to z_0}f(z)=\lim_{\epsilon\to\epsilon_0}\oint f(z)ds/\oint ds$.

When $s$ is the length of a circular path of radius $\epsilon$ surrounding $z_0$:

$s = \epsilon\theta$.

$\oint f(z)ds=\int_0^{\theta_{\max}}f(z_0+\epsilon e^{i\theta})\epsilon d\theta$.

$\oint ds=\theta_{\max}\epsilon$.

$\lim_{z\to z_0}f(z)=\lim_{\epsilon\to\epsilon_0}\frac{1}{\theta_{\max}}\int_0^{\theta_{\max}}f(z_0+\epsilon e^{i\theta})d\theta$\\
\\
To see how it works, let $I=\int_{0}^{\theta_{\max}}f(z)d\theta$ and consider $f(z) = z^\alpha$.

When $z_0 \neq 0$,

$f(z)=(z_0+\epsilon e^{i\theta})^\alpha=z_0^{\alpha}+\alpha z_0^{\alpha-1}\epsilon e^{i\theta}+O(\epsilon^2)$

$I=z_0^{\alpha}\theta_{\max}+\epsilon\int_0^{\theta_{\max}}\alpha z_0^{\alpha-1}e^{i\theta}d\theta+O(\epsilon^2)$

$\lim_{z\to z_0}f(z)=\lim_{\epsilon\to\epsilon_0}I/\theta_{\max}$

so the leading term limit when $z_0 \neq 0$ is

$^1\lim_{z\to z_0}z^{\alpha}=z_0^{\alpha}$\\
\\
When $z_0 = 0$ and $\alpha \neq 0$,

$f(z)=\epsilon^{\alpha}e^{i\alpha\theta}$

$I=\frac{\epsilon^{\alpha}}{i\alpha}\left [ e^{i\alpha\theta}\right]_0^{\theta_{\max}}$

As with reals, take the mean of an oscillating function over a period.

$\theta_{\max} = 2\pi/\alpha$
 
$\lim_{z\to 0}z^{\alpha}=0$\\
\\
When $\alpha = 0$,
$\lim_{z\to z_0}1=1$\\
\\
Defining $s$ as path length differs from normal contour integration. It is important not to get the two mixed up.

$\oint f(z)ds=\int_0^{\theta_{\max}}f(z_0+\epsilon e^{i\theta})\epsilon d\theta$.

$\oint f(z)dz=\int_0^{\theta_{\max}}f(z_0+\epsilon e^{i\theta})(i\epsilon e^{i\theta})d\theta$.

\section{Improper Integrals}

Evaluating definite integrals using limits of series requires care and whenever feasible the problem that yields the integral should be examined for clues to the behaviour near infinity and near each singularity.

When no such information is available and when an indefinite integral is available, substitute $\omega$ for infinity and when the limit is at a finite singularity $x_0$ substitute $x_0 \pm 1/\omega$.

A Riemann sum can often be used to calculate the leading term limit of a definite integral. Terms beyond the leading term limit should sometimes be treated with suspicion and where high accuracy is needed a cubic spline fit~\cite{31} would be better. The following Riemann sum can be used for any function without singularities on domain $[0, \omega)$.
\begin{eqnarray*}
&\int_0^{\omega}f(x)dx=\lim_{m\to\omega}\left(\frac{1}{m}\lim_{n\to m\omega}\sum_{k=1}^nf((k-.5)/m)\right)\\
&=\lim_{n\to\omega^2}\frac{1}{\sqrt{n}}\sum_{k=1}^nf(x)\mbox{  when  }x=(k-.5)/\sqrt{n}
\end{eqnarray*}
For example, using the new definition of ``lim", this immediately yields the exact results:
\begin{eqnarray*}
&\int_0^{\omega}dx=\omega\qquad\int_0^{\omega}xdx=\omega^2/2\qquad\int_0^{\omega}\cos(x)dx=0\\
&\int_0^{\omega}\sin(x)-x\cos(x)dx=0
\end{eqnarray*}
More commonly it yields approximate results:
\begin{eqnarray*}
& \int_0^{\omega}\exp(x)dx=\frac{\exp(1/2\omega)(\exp(\omega)-1)}{\omega(\exp(1/\omega)-1)}\approx\frac{(1+1/2\omega+1/8\omega^2)(\exp(\omega)-1)}{(1+1/2\omega+1/6\omega^2)}
& \approx\exp(\omega)
\end{eqnarray*}
This has the correct leading term limit, with an error term smaller than the leading term by a factor of $O(\omega^2)$. That's a typical order of magnitude factor for an error term of a Riemann sum.

When there is a singularity on the reals, a transformation can help. Let $x = t(\xi)$ where $\xi$ is divided into equal intervals and each interval is evaluated at its midpoint. The width of each $x$ interval is $\Delta x =t'(\xi)\Delta\xi$. Then
\begin{displaymath}
\int_{x_{\min}}^{x_{\max}}f(x)dx=\lim_{n\to\omega^2}\frac{1}{\sqrt{n}}\sum_{k=k_{\min}}^{k_{\max}}f(t(\xi))t'(\xi)
\end{displaymath}
where $\xi=(k+k_{\min}-0.5)/\sqrt{n}$, and $x_{\min} = t(k_{\min})$ and $x_{\max} = t(k_{\max})$. Change the sign if limits are swapped.

Useful transformations include:
\begin{itemize}
\item	$t(\xi) = x_{\min} + 1/\xi$.
This maps $(0, \omega)$ to $(x_{\min}, \omega)$ for functions that misbehave at $x_{\min}$ but integrate to a finite number at the infinite limit.
\item	$t(\xi) = x_{\min} + \xi /(1-\exp(-\xi))$
This maps $(-\omega, \omega)$ to $(x_{\min}, \omega)$ for functions that misbehave at $x_{\min}$ and infinity.
\item	$t(\xi) = (x_{\max} + x_{\min}) / 2 + (x_{\max} - x_{\min}) \tan^{-1}(\pi\xi / (x_{\max}-x_{\min}) ) / \pi$.
This maps $(-\omega, \omega)$ to $(x_{\min}, x_{\max})$ for functions that misbehave at both limits.
\end{itemize}
\section{Dirac delta function}

I don't want to deal in any detail with generalized functions in this paper, but there's no difficulty in squaring the Dirac delta function~\cite{32}. For generalized functions that can be constructed from limits, everything remains commutative and there's no need for non-commutative algebras.\\
\\
Let $f(x)$ be any probability density function on the real numbers. Let $\epsilon_1$ be any appropriate infinitesimal number. Then the Dirac delta function is the set of functions defined by $\delta(x)=\lim_{\epsilon\to\epsilon_1}f(x/\epsilon)/\epsilon$. To see that the limit exists, express it as the limit of a smooth sequence as follows. Let $g(\omega)=-1/\epsilon_1$, then the limit becomes $\lim_{n\to\omega}g(n)f(xg(n))$. Call this limit $f(x/\epsilon_1)/\epsilon_1$. The square of the Dirac delta function is simply $f(x/\epsilon_1)^2/\epsilon_1^2$. To get the correct value of the Heaviside function at zero, $H(0)=0.5$, the initial probability density function should be symmetric about zero.

\section{Conclusions}

I hope that the methods of this paper go some way towards banishing the concept of divergence. Methods are presented here for evaluating the limit of sequences that converge/diverge to infinity, for evaluating limits of oscillatory and general sequences, for finding limits on the reals, finding limits on complex numbers, and for evaluation of improper integrals.

Without exception, each of the methods used here is ``local" on the integers and real numbers in the sense that the mean is defined at a specific integer or real number without recourse to adjacent integers or real numbers. The ensemble mean of a sequence at given index $n$ is totally independent of the values of the sequence at values $n-1$ and $n+1$. This takes some getting used to, and is possible only by using the properties of the mapping that is used to generate the sequence. On complex numbers this restriction is relaxed to defining the mean using an integral over a circle of fixed infinitesimal ``local" radius.

This part of this paper leaves several questions unanswered. Although proofs aren't needed for practical application of the methods presented here, ideally proper proofs/disproofs should be found for the conjectures:
\begin{itemize}
\item	All infinite sequences can be decomposed into a finite sum of smooth, periodic, chaotic and random components and this decomposition yields a unique leading term limit.
\item	Every definite integral on the reals has a unique value on $\mathbb{I}_n$.
\end{itemize}

As a piece of applied mathematics, the number one priority would be in codifying the methods presented here and extensions of them into software, so that someone could ask, for example ``$\int_0^{\infty}\sin(x)dx=$ ?", and get back an immediate exact answer.


\begin{thebibliography}{32}
\bibitem{Part1} D.A. Paterson (2011) Banishing divergence Part 1: Infinite numbers as the limit of sequences of real numbers. arXiv submission 0304813.
\bibitem{bijection} A. Bertiger (2009) Different Sizes of Infinity Lecture 2: Bijections http: //www.math.cornell.edu/$\sim$mec/2008-2009/Bertiger/Lecture\_2.html
\bibitem{BigO} http://en.wikipedia.org/wiki/Big O notation
\bibitem{2}  K.D. Stroyan and W.A.J. Luxemburg (1976), Introduction to the theory of infinitesimals, Academic Press, New York.
\bibitem{25}  http://en.wikipedia.org/wiki/Ramanujan\_summation
\bibitem{26}  G. Box and G. Jenkins (1970) Time series analysis: Forecasting and control, Holden-Day, San Francisco.
\bibitem{27}  http://en.wikipedia.org/wiki/Stationary\_ergodic\_process
\bibitem{ensemble} R. Sasaki (2010), http://www.stanford.edu/$\sim$rsasaki/EEAP248/slide1
\bibitem{28}  http://en.wikipedia.org/wiki/Ces\%C3\%A0ro\_summation
\bibitem{29}  http://en.wikipedia.org/wiki/H\%C3\%A9non\_map
\bibitem{lebesgue} http://planetmath.org/encyclopedia/Integral2.html
\bibitem{30}  http://mathworld.wolfram.com/PrimeNumberTheorem.html
\bibitem{orderstats} http://mathworld.wolfram.com/OrderStatistic.html
\bibitem{31}  http://www.physics.utah.edu/$\sim$detar/phys6720/handouts/cubic\_spline/
cubic\_spline/node1.html
\bibitem{32}  http://en.wikipedia.org/wiki/Dirac\_delta\_function

\end{thebibliography}
\end{document}